\documentclass[12pt]{article}
\usepackage{latexsym}
\usepackage{amsmath, amsthm}
\usepackage{graphicx}
\usepackage{amssymb}

\newtheorem{thm}{Theorem}[section]

\newtheorem{prop}[thm]{Proposition}

\title
{Two analogs of intrinsically linked graphs}

\author{Chris Cicotta \and Joel Foisy  \and Tom Reilly  \and Sara Rezvi  \and Ben Wang  \and Alice Wilson}

\begin{document}

\maketitle

\begin{abstract}
A graph $G$ is {\it intrinsically $S^1$-linked} if for every embedding of the vertices of $G$ into $S^1$, vertices that form the endpoints of two disjoint edges in $G$ form a non-split link in the embedding.  We show that a graph is intrinsically $S^1-$linked if and only if it is not outer-planar.  A graph is {\it outer-flat} if it can be embedded in the $3-$ball such that all of its vertices map to the boundary of the $3-$ball, all edges to the interior, and every cycle bounds a disk in the $3-$ball that meets the graph only along its boundary.  We show that a graph is outer-flat if and only if it is planar. 


\end{abstract}
\section{Introduction}

Recall that a graph is said to be {\it intrinsically linked} if it contains a pair of non-splittably linked cycles in every (tame) spatial embedding.  Conway, Gordon \cite{Conway} and Sachs \cite {Sachs1} showed that the graph $K_6$ is intrinsically linked.  Sachs went on to conjecture that the Petersen Family of graphs, the $7$ graphs obtained from $K_6$ by triangle-Y and Y-triangle exchanges, form the complete set of minor-minimal intrinsically linked graphs.  Robertson, Seymour and Thomas later proved Sachs' conjecture (see \cite {Robertson} for an outline of the proof and other references).  They further showed that a graph has a linkless embedding if and only if it has a flat embedding (every cycle of $G$ bounds a disk that meets the graph only along its boundary).  In this paper, we define two analogous versions of intrinsically linked.  The first is a lower-dimensional analog called intrinsically $S^1$-linked.  We say that a graph $G$ is {\it intrinsically $S^1$-linked} if for every embedding of the vertices of $G$ into $S^1$, vertices that form the endpoints of two disjoint edges in $G$ form a non-split link in the embedding.  We later show that a graph is intrinsically $S^1$-linked if and only if it is not outerplanar.  The other version of intrinsically linked, we will call intrinsically outer-linked.  We will define this concept later in the introduction.  Closely related to this is the concept of outer-flat.  A graph is {\it outer-flat} if it can be embedded in the $3-$ball such that all of its vertices map to the boundary of the $3-$ball, all edges to the interior, and every cycle bounds a disk in the $3-$ball that meets the graph only along its boundary.  We show that a graph is outer-flat if and only if it is planar. 

These analogues of intrinsically linked were motived in part by some characteristics of the Colin de Verdi\`ere graph parameter, $\mu(G)$.  This parameter, though defined in terms of spectral properties of a graph, turns out to be closely related to topological properties of a graph (see \cite {holst} for a survey on this parameter).  In particular, by \cite{dev1} and \cite{dev2}, a graph is outerplanar if and only if $\mu(G) \leq 2$, planar if and only if $\mu(G) \leq 3$, and, by \cite{lovasz}, possesses a flat embedding if and only if $\mu(G) \leq 4$.

Here we establish some terminology relevant for the definition of intrinsically $S^1$-linked.  Recall that a link in $S^1$ is an embedding of two copies of $S^0$ into $S^1$. Let $\{p_1, p_2\}$ and  $\{p_3, p_4\}$ denote the two copies of $S^0$ that make up a link in $S^1$.  We say $\{p_1, p_2\}$ and  $\{p_3, p_4\}$ form a {\it non-splittable link in $S^1$} provided $p_1$ and $p_2$ lie in different components of  $S^1$ - \{$p_1$, $p_2$\}, otherwise we say the link is {\it splittable}.  An {\it $S^1$-embedding} of a graph $G$ is an injective map of all the vertices of G into $S^1$.  Given an $S^{1}$-embedding of a graph $G$, we will consider two vertices to form an $S^0$ if and only if they are the end points of an edge.  Finally, we may say that a graph $G$ is intrinsically $S^1-$linked if every $S^1$-embedding of $G$ contains at least one pair of disjoint copies of $S^0$ that form a non-splittable link.

Recall that a graph $G$ is {\it outerplanar} if it can be embedded in the plane so that all its vertices touch the outer face, and recall that a graph $H$ is said to be 
a {\it minor} of a graph $G$ if $H$ can be obtained from $G$ by 
deleting and/or contracting a finite number of edges of $G$.   Halin \cite {halin} (see also Chartrand and Harary \cite{chartrand}) first established the well-known result that a graph is outerplanar if and only if it contains $K_4$ or $K_{3,2}$ as a minor.  We use this result to show that a graph is outerplanar if and only if it is not intrinsically $S^1$-linked.

In this last part of the introduction, we establish some terminology to lead up to the definition of intrinsically outer-linked.  Recall that an alternate definition for outerplanar is for a graph $G$ to have an embedding into the disk such that all vertices of the graph lie on the boundary of the disk.  We analogously define an {\it outer-embedding} of a graph $G$ to be an embedding  of $G$ into the $3-$ball such that all of the vertices of $G$ map to the boundary of the $3-$ball, and all edges map to the interior.  In an outer-embedding of a graph, a {\it splittable link in $B^{3}$} consists of a cycle, $\gamma$, and an edge, $e$, for which there exists a disk that is homeomorphically mapped to $B^3$, such that the boundary of the disk maps to the boundary of $B^3$, and the interior to the interior, and the disk cuts $B^3$ into two components, one containing $e$, and the other $\gamma$.  We call a link {\it non-splittable} (or just non-split) otherwise.  Finally, we can define a graph $G$ to be {\it intrinsically outerlinked} if every outer-embedding of $G$ contains a non-split link consisting of a cycle and an edge.  Our second main result in the paper will be that a graph is intrinsically outer-linked if and only if it is not planar.  We call a graph {\it outer-linkless} provided it has an outer-embedding for which every link in $B^3$ is splittable.   We further show that a graph is outer-flat if and only if it is outer-linkless.



In a future work, we will examine the problem of classifying all graphs that have a disjoint pair of non-split links in every $S^1$-embedding, and the problem of classifying all graphs that have a disjoint pair of non-split outer-links in every outer-embedding.  The analogous problem of determining the complete set of minor-minimal graphs that contain a disjoint pair of non-splittable links in every spatial embedding was studied in \cite{chan}.  Though the set of such graphs is finite (due to the powerful result in \cite {robertson1}), and some minor-minimal examples were given in \cite{chan}, the problem seems very difficult at this time.  It is our hope that the analogous problems in the contexts of $S^1$-embeddings and outer-embeddings will be more tractable, and will shed some light on the spatial problem.

\section {Intrinsically $S^{1}$-linked Graphs}

\begin{thm}
The graphs $K_4$ and {$K_{3,2}$ are intrinsically $S^{1}$-linked.}
\end{thm}
\begin{proof} 
Up to symmetry, there is only one $S^1$-embedding $K_4$, and in that embedding, there is a non-split link.  Up to symmetry, there are two ways to place the five vertices of $K_{3,2}$ on $S^{1}$. In both cases, there is a non-split link. Therefore, $K_{3,2}$ is intrinsically $S^{1}$-linked.

\end{proof}

Now we present an alternative proof of  intrinsic $S^{1}$-linking of $K_{3,2}$ using an analog of Conway-Gordon and Sach's proof that $K_6$ is intrinsically linked.  An analogous proof can also be constructed for $K_4$, which we will omit.  Indeed, such a proof was given for $K_4$ in \cite{taniyama} as part of a more general theorem.  For our proof for $K_{3,2}$, we need to define the mod $2$ linking number for links in $S^1$.  The mod $2$ linking number of two copies of $S^{0}$ embedded in $S^1$ is $1$ if they are non-splittably linked and a $0$ otherwise.


%

%


\begin{proof}  

We will denote the vertices of $K_{3,2}$ as $\{a,b,c,1,2\}$, where $\{a,b,c\}$ make up one partition.  Consider the $S^1$-embedding of $K_{3,2}$ that places the vertices in the clock-wise order $a,b,c,1,2$.  For this $S^1$-embedding, the sum of the mod $2$ linking numbers is odd.  One can obtain any other $S^1$-embedding by a sequence of exchanging vertices pair-wise.  It thus suffices to show that a single exchange of vertices leaves the sum of mod $2$ linking numbers odd.  For the one exchange, we can switch vertices within a partition, or switch vertices in different partitions.  Clearly, exchanging vertices within the same partition will not affect the sum of the mod $2$ linking numbers on all pairs of links.  So we consider the affect of switching two vertices in different partitions.  Without loss of generality, we consider the affect of switching the vertices $a$ and $1$.  The only links that will be affected use both $a$ and $1$.  These are exactly two pairs:  $(a2,1b)$ and $(a2,1c)$.  Each pair will experience a change of plus or minus one in linking number, thus the total change in the sums of all of the mod $2$ linking numbers will be even.  It follows that the sum of the mod $2$ linking numbers will remain odd.  Therefore, $K_{3,2}$ is intrinsically $S^{1}$-linked.
\end{proof}

Finally, we have:

\begin{prop} If a graph $G$ is outerplanar then $G$ is not intrinsically $S^{1}$-linked.
\end{prop}

\begin{proof}
If a graph $G$ is outerplanar, then we embed $G$ in a disk, with all of the vertices of $G$ on the boundary of the disk.  Restricting this embedding to the vertices of $G$ gives a outer-linkless $S^1$-embedding of $G$.
\end{proof}

For the other direction, we need the following well-known result:

\begin{thm} \cite{halin}  A graph is outerplanar if and only if it does not contain $K_{4}$ nor $K_{3,2}$ as a minor. 
\end{thm}

We will have shown that a graph is outerplanar if and only if it is not intrinsically $S^1$-linked after we have established the following theorem. 

\begin{thm} Vertex expansion preserves intrinsic $S^{1}$-linking.  
\end{thm}

\begin{proof}

Consider a graph $G$ that is intrinsically $S^{1}$-linked. Consider the graph $G' $, obtained from $G$ by a vertex expansion, where vertex $v$ is expanded to edge $(v', v'')$.  We claim that $G'$ is intrinsically $S^{1}$-linked.

Consider an arbitrary $S^{1}$-embedding of $G'$.  We will show that it contains a non-split link. We consider the associated $S^{1}$-embedding of $G$ that results from omitting the vertex $v''$ from the $S^1$-embedding of $G'$, and relabeling $v'$ as $v$.  In the first case, suppose $v$ does not form part of a non-split link.  In this case, placing $v''$ back on the circle will not destroy any non-split links, and thus the $S^1$-embedding of $G'$ contains a non-split link. 

In the second case, we must consider when $v$ is involved in a non-split link.  Let us say the vertices of the edges $(v,w)$ and $(a,b)$ represent the two copies of $S^0$ that make up this link.  Let us consider what happens after we put $v''$ back on the circle to get back to the original $S^1$-embedding of $G'$.  There are two cases two consider.  In the  first case, $v''$ and $w$ lie in the same component of $S^1-\{a,b\}$.  In this case, a non-split link in the $S^1$-embedding of $G'$ is $(v',v'')$ and $(a,b)$.  In the second case, $v''$ and $w$ lie in different components of $S^1-\{a,b\}$.  In this case, the non-split link in the $S^1$-embedding of $G'$ is either $(v',w)$ or $(v'',w)$, depending on which of $v'$ and $v''$ is adjacent to $w$ in $G'$.

Thus, in all cases, we can find a non-split link in the $S^1$-embedding of $G'$. Therefore vertex expansion preserves intrinsic $S^{1}$-linking.

%

\end{proof} 

Since a non-outerplanar graph can be obtained from $K_4$ or $K_{3,2}$ by a sequence of vertex expansions and/or adding vertices and/or edges (and the latter two operations clearly preserve intrinsic $S^1$-linking), we have thus established our main theorem for this section:

\begin{thm}  A graph is intrinsically $S^1$-linked if and only if it is not outerplanar.
\end{thm}

  \section{Intrinsically outerlinked graphs}

\begin{prop} If a graph $G$ is planar, then $G$ is outer-flat and outer-linkless.
\label{flat}
\end{prop}
\begin{proof}
 Let $G$ be a planar graph.  Embed $G$ on the boundary of $B^3$.  Deform the edges slightly so that they lie in the interior of $B^3$.  The resulting outer-embedding  of $G$ is outer-flat and outer-linkless.
 
 \end{proof}
 
\begin {figure}
 
 \vskip -2.2in
 
 \includegraphics [scale=.5]{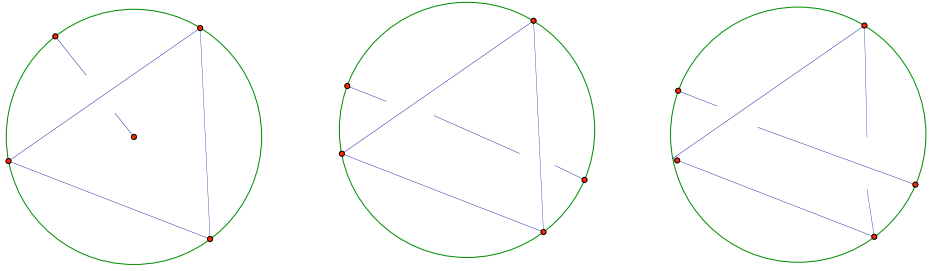}
 
 \vskip -1.3in
 \caption{A non-calculable projection, and two calculable projections.  In the second projection, the mod $2$ linking number is $0$.  In the last projection, the mod $2$ linking number is $1$.}
 \label{calc}
 \end{figure}

 In our study of intrinsically outer-linked graphs, a useful tool for us will be an analog of the mod $2$ linking number.  In an outer-embedding of a graph, a {\it link in $B^{3}$} consists of a cycle, $\gamma$, and a disjoint edge, $e$.  We allow deformations of links that are analogous to deforming links in space, except that during the deformations, vertices must stay on the boundary of $B^3$.  We say that $\gamma$ and $e$ have non-zero mod $2$ linking number if the image of $\gamma$ can be extended to the image of a disk with an odd number of transversal intersections with $e$.  In order to calculate mod $2$ linking number, not only do we need a regular projection, but we need the further restriction that the projection be {\it calculable}.  We say that a projection of a link in $B^3$ is calculable provided that the vertices of the edge lie on the equator of the boundary of $B^3$ (otherwise there is an ambiguity as to whether or not the vertex lies in the upper or lower hemisphere, see Figure \ref {calc}).  Given a calculable projection, we compute the mod $2$ linking number of $\gamma$ and $e$ to be the number of times, mod $2$ that $e$ crosses over $\gamma$ in that projection.  If the mod $2$ linking number of $\gamma$ and $e$ is non-zero, then $\gamma $ and $e$ form a non-splittable link in $B^3$.
 
\begin{thm} The graphs $K_{5}$ and $K_{3,3}$ are intrinsically outerlinked.
\label{label}
\end{thm}
\begin{proof}
We present two proofs.  First, consider an arbitrary outer-embedding of $K_{5}$.  Deform the vertices of $K_{5}$ to the equator of the sphere, and then embed the $B^3$ into space.  We can place a new sixth vertex outside the ball, connecting it to the first five vertices, and thus make an embedding of $K_{6}$.  We can connect the sixth vertex in such a way that in a projection of our embedding of $K_6$, the new edges have no crossings with the edges from the original $K_5$.  Since $K_6$ is intrinsically linked, and contains a pair of cycles linked with non-zero mod $2$ linking number \cite {Conway} \cite{Sachs1}, there are a pair of linked cycles in the $K_6$ embedding, say $C_1$ and $C_2$, where $C_2$ consists of edges $e_1,e_2,$ and $e_3$, where $e_2$ and $e_3$ are incident to the sixth vertex.  Since the edges from the sixth vertex have no crossings, all of the crossings must take place inside the projection of the $B^3$, thus $C_1$ and $e_1$ formed a non-split link (with mod $2$ linking number $1$) in the outer-embedding of $K_5$.  Thus $K_5$ is intrinsically outer-linked.  Since $K_{3,3,1}$ is intrinsically linked with a pair of cycles with non-zero mod $2$ linking number in every spatial embedding \cite {Sachs1}, a similar proof works for $K_{3,3}$.
\end{proof}

\begin{figure}[h]
\begin{center}
\includegraphics[scale=.2]{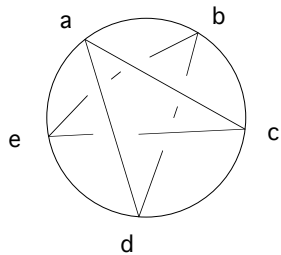}
\vskip -.25in
\caption{A projection of an outer-embedding of $K_5$.}
\label{k5}
\end{center}
\end{figure}

\begin{proof}
Now we present a proof of the Theorem that is analogous to Conway-Gordon and Sachs' proof that $K_6$ is intrinsically linked.  Consider the outer-embedding of $K_5$ represented in Figure \ref{k5}.    Of all of the $10$ possible links, the only non-split link is $(abd,ec)$, with mod $2$ linking number equal $1$.  Thus, the sum of mod $2$ linking numbers of all of the links is $1$.  To go from one outer-embedding to another, we allow for the usual types of deformations (with vertices allowed to move around only on the boundary), but we must also allow for crossing changes.  We claim that a change of crossings will leave invariant the mod $2$ sum of mod $2$ linking numbers of all possible links.  The non-trivial case consists of a crossing change between non-adjacent edges.  In this case, there is exactly one vertex that is not incident to either edge.  This vertex, and the edges, form two possible links consisting of a cycle and an edge.  Each cycle and edge pair will have its mod $2$ linking number change by exactly plus or minus $1$.  Thus, the sum of the mod $2$ linking numbers will change by either plus or minus $2$, or $0$.  Thus the mod $2$ sum will be $1$ for every outer-embedding of $K_5$.  The proof for $K_{3,3}$ is similar.
\end{proof}

\vskip -.45in

\begin{thm}  If a graph $G$ is non-planar, then $G$ is intrinsically outer-linked.
\label{non-planar}
\end{thm}

\begin{proof} 
If $G$ is non-planar, then $G$ contains $K_5$ or $K_{3,3}$ as a minor (\cite{kuratowski}, \cite{wagner}).  Let $*$ denote the usual join of a graph and a vertex.  As in the first proof of Theorem \ref{label}, an outer-embedding of $G$ induces a spatial embedding of $G*v$, which has a regular projection for which there are no crossings on edges incident to $v$.  Moreover, the graph $G*v$ contains $K_6$ or $K_{3,3,1}$ as a minor, and thus contains a pair of non-splittably linked cycles, with non-zero mod $2$ linking number (\cite {fellows}, \cite{nes}).  Again, as in the first proof of Theorem \ref{label}, there must exist a non-splittable link in the outer-embedding of $G$.
\end{proof}

 \noindent
 We have established our main theorem of this section:
 
\begin{thm} The following are equivalent for a graph $G$:\newline
\indent \indent 1. $G$ is planar\newline
\indent \indent 2. $G$ is outer-flat.\newline
\indent \indent 3. $G$ is outer-linkeless.\newline
\end{thm}

  \noindent
{ \bf  Acknowledgment}
  
The results in this paper were obtained by a research group in the 2005 Research Experience for Undergraduates at SUNY Potsdam, advised by Joel Foisy and sponsored by National Science Foundation Grant DMS-0353050 and National Security Grant MDA H982300510095.

  	\bigskip
	
		{\footnotesize  CHRISTOPHER CICOTTA, CLARKSON UNIVERSITY, POTSDAM, NY 13699}
		{\footnotesize {\it E-mail address:}  cicottcm@clarkson.edu}

		{\footnotesize  JOEL FOISY, DEPARTMENT OF MATHEMATICS, SUNY POTSDAM, 
		POTSDAM, NY 13676}
		{\footnotesize {\it E-mail address:}  foisyjs@potsdam.edu}
		
		{\footnotesize  SARA REVZI, UNIVERSITY OF CHICAGO, CHICAGO, IL 60637}
		{\footnotesize {\it E-mail address:}  arsinoe@uchicago.edu}
		
		{\footnotesize  TOM REILLY, UNION COLLEGE, SCHENECTADY, NY 12308}
		{\footnotesize {\it E-mail address:}  reillyt@union.edu}
		
		{\footnotesize  BEN WANG, BINGHAMTON UNIVERSITY, BINGHAMTON, NY 13902}
		{\footnotesize {\it E-mail address:}  benactuarial@gmail.com}
		
		{\footnotesize  ALICE WILSON, SUNY POTSDAM, POTSDAM, NY 13676}
		{\footnotesize {\it E-mail address:}  wilson09potsdam.edu}

\end{document}